\documentclass[]{article}
\begin{document}
\newtheorem{proposition}{Proposition}[section]
\newtheorem{definition}{Definition}[section]
\newtheorem{lemma}{Lemma}[section]

\title{\bf Zero Divisor Manifolds}
\author { Keqin Liu\\Department of Mathematics\\The University of British Columbia\\Vancouver, BC\\
Canada, V6T 1Z2}
\date{March, 2025}
\maketitle

\begin{abstract} We develop the basic properties of 
$R^{(2)}$-modules,  introduce  the concept of zero divisor manifolds, 
construct projective  $R^{(2)}$-space which generalizes the real projective space, and initiate the study of the counterpart of symplectic spaces.
\end{abstract}

\bigskip
The theory of real manifolds has some natural generalizations if Euclidean spaces are replaced by the modules over some real commutative associative algebras with zero divisors. The simplest generalization of real manifolds obtained in this way  is zero divisor manifolds introduced in this paper. Instead of looking like a 
Euclidean space, a  zero divisor manifold looks-like locally a module over the dual real number algebra $R^{(2)}$, where $R^{(2)}$ is a real commutative  associative algebras with zero divisors and was introduced by W. K. Clliford in 1873.  After developing the basic properties of 
$R^{(2)}$-modules, we introduce  the concept of zero divisor manifolds, 
construct projective  $R^{(2)}$-space which generalizes the real projective space, and initiate the study of the counterpart of symplectic spaces.

\bigskip
\section{Modules over the Dual Real Number Algebra}

 Let $\mathcal{R}$ be the real number field with the identity $1$. The dual real number algebra is the 2-dimensional real associative algebra
$\mathcal{R}^{(2)}=\mathcal{R}\oplus \mathcal{R}\,1^{\#}$ with 
a basis $\{1, \,1^{\#}\}$ and  the associative product defined by 
$$
(x_1+x_2\,1^{\#})(y_1+y_2\,1^{\#})=x_1y_1+(x_1y_2+x_2y_1)\,1^{\#},
$$
where $x_1$, $x_2$, $y_1$, $y_2$ are real numbers. If 
$x=x_1+x_2\,1^{\#}\in \mathcal{R}^{(2)}$ with $x_1$, $x_2\in \mathcal{R}$,
then $Re\,x:=x_1$ and $Ze\,x:=x_2$ are called the {\bf  real part} and the  {\bf zero-divisor part} of $x$, respectively. The dual real number algebra $\mathcal{R}^{(2)}$ has many zero-divisors.
In fact, if $0\ne x\in \mathcal{R}^{(2)}$, then $x$ is a zero-divisor if and only if $Re\,x=0$,  and
$x$ is invertible if and only if $Re\,x\ne 0$. Moreover if $x$ is invertible,  then  the inverse $\displaystyle\frac{1}{x}:=x^{-1}$ of $x$ is given by
$\displaystyle\frac{1}{x}=x^{-1}=\displaystyle\frac{1}{Re\, x}-\displaystyle\frac{Ze\, x}{(Re\, x)^2}\,1^{\#}$. For convenience,  the product $x^{-1}y$ is also denoted by $\displaystyle\frac{y}{x}$, where $x$, $y\in \mathcal{R}^{(2)}$ and  $x$ is invertible.
The  dual real number algebra $\mathcal{R}^{(2)}$ is a normed algebra with respect to the following norm:
$$
||x||:=\sqrt{2\,(Re\,x)^2+(Ze\,x)^2}\quad\mbox{for $x\in\mathcal{R}^{(2)}$.}
$$

\medskip
A  $\mathcal{R}^{(2)}$-module $V$ has the following submodule chain:
$$\{0\}\subseteq Im\,1^{\#}_V \subseteq Ker\,1^{\#}_V\subseteq V,$$
where the submodules $Im\,1^{\#}_V $ and $ Ker\,1^{\#}_V$ are defined by
$$
Im\,1^{\#}_V:=\{ 1^{\#}\,v\,|\, v\in V\}, \qquad 
Ker\,1^{\#}_V:=\{ v\,|\, \mbox{$1^{\#}\,v=0$ and $v\in V$}\}.
$$
An element of a $\mathcal{R}^{(2)}$-module $V$ is  called a {\bf dual real vector} of $V$. 

\medskip
We now give the following generalization of the concept of bases in Euclidean space.

\medskip
\begin{definition}\label{def1} A subset $\beta$ of  a 
$\mathcal{R}^{(2)}$-module $V$ is called a 
{\bf $(\mathcal{R}^{(2)}, \mathcal{R})$-basis} for $V$ if $\beta$ has the following three properties:
\begin{description}
\item[(i)] $\beta=(S_1\,\|\, S_2):=S_1\cup S_2$, where $S_1\subseteq V$  and $S_2\subseteq   Ker\,1^{\#}_V$;
\item[(ii)]  $\beta=(S_1\,\|\, S_2)$ is  {\bf $(\mathcal{R}^{(2)}\,\|\, \mathcal{R})$-linearly independent}, i.e., for any 
finite subset $(\{v_i\}_{i=1}^s\,\|\, \{v_{n+j}\}_{j=1}^t)$ of $\beta$ with $\{v_i\}_{i=1}^s\subseteq S_1$ and 
$\{v_{n+j}\}_{j=1}^t\subseteq S_2$, we have
\begin{eqnarray*}
&&\sum_{k=1}^{s+t}\,c_kv_k=0,\quad  \{c_i\}_{i=1}^{s}\subseteq \mathcal{R}^{(2)} \quad\mbox{and}\quad
\{c_{n+j}\}_{j=1}^{t}\subseteq \mathcal{R}\\
&&\Longrightarrow c_k=0 \quad\mbox{for $1\le k\le s+t$};
\end{eqnarray*}
\item[(iii)] For any dual real vector $v$ of $V$, there exist  a finite number of  dual real numbers 
$\{c_k\}_{k=1}^{s+t}$ and a finite number of  dual real vectors
$\{v_k\}_{k=1}^{s+t}$  such that  
$$
\{c_i\}_{i=1}^{s}\subseteq \mathcal{R}^{(2)},  \quad
\{c_{n+j}\}_{j=1}^{t}\subseteq \mathcal{R}\quad\mbox{ and}\quad  v=\sum_{k=1}^{s+t}\,c_kv_k.
$$
\end{description}
\end{definition}

\medskip
\begin{proposition}\label{pr1} Every $\mathcal{R}^{(2)}$-module has a 
$(\mathcal{R}^{(2)}, \mathcal{R})$-basis.
\end{proposition}

\medskip 
If a $\mathcal{R}^{(2)}$-module $V$ has a finite $(\mathcal{R}^{(2)}, \mathcal{R})$-basis
$(\{v_i\}_{i=1}^n\,\|\,\{v_{n+j}\}_{j=1}^m)$, then $V$ is called {\bf finite dimensional}, the 
unique pair $(n, m)$ of non-negative integers is called the $\mathcal{R}^{(2)}$-{\bf dimension} of $V$ and is denoted by 
$\dim_{\mathcal{R}^{(2)}}(V)=(n,\,m)$. 

\medskip 
Let  $n$ and $m$ be nonnegative integers with $n+m\ne 0$. Let
$$\mathcal{R}^{(2)n,m}:=
\left\{(x^1,  \dots ,  x^n\,| \,x^{n+1}1^{\#},  \dots ,  x^{n+m}1^{\#})\left |
\begin{array}{c}
x^1, \dots, x^n\in \mathcal{R}^{(2)}\\ 
x^{n+1},  \dots ,  x^{n+m}\in \mathcal{R}\end{array}\right.\right\}.$$
Then $\mathcal{R}^{(2)n,m}$  is a  $\mathcal{R}^{(2)}$-module with the operation of coordinatewise addition and scalar multiplication. An element in $\mathcal{R}^{(2)n,m}$ is called a {\bf real $(n, m)$-vector}, and  
$\mathcal{R}^{(2)n,m}$ is called 
the {\bf real $(n, m)$-dimensional $\mathcal{R}^{(2)}$-module}. Let
\begin{eqnarray*}
e_1:&=&(\underbrace{1, 0, \dots ,  0, 0}_{n}\,| \,\underbrace{0, 0, \dots , 0, 0}_{m})\\
&\vdots&\\
e_n:&=&(\underbrace{0, 0, \dots ,  0, 1}_{n}\,| \,\underbrace{0, 0, \dots , 0, 0}_{m})\\
e_{n+1}:&=&( \underbrace{0, 0\dots ,  0, 0}_{n}\,| \,\underbrace{1^{\#}, 0, \dots , 0, 0}_{m})\\
&\vdots&\\
e_{n+m}:&=&( \underbrace{0, 0\dots ,  0, 0}_{n}\,| \,\underbrace{0, 0, \dots , 0, 1^{\#}}_{m}).
\end{eqnarray*}
Then $(\{e_i\}_{i=1}^n\,\|\, \{e_{n+j}\}_{j=1}^m)$ is a 
$(\mathcal{R}^{(2)}, \mathcal{R})$-basis, which is called  the {\bf standard $(\mathcal{R}^{(2)}, \mathcal{R})$-basis} for the  real $(n, m)$-module $\mathcal{R}^{(2)n, m}$.

\medskip
\begin{proposition}\label{pr1.2} Let $V$ be a $\mathcal{R}^{(2)}$-module. Then $V$ is isomorphic to
 $\mathcal{R}^{(2)n, m}$ if and only if $\dim_{\mathcal{R}^{(2)}}(V)=(n,\,m)$.
\end{proposition}

\bigskip
\section{$\mathcal{R}^{(2)}$-Differentiability}

Let $(\{e_i\}_{i=1}^n || \{e_{n+j}\}_{j=1}^m)$ be  the standard 
$(\mathcal{R}^{(2)}, \mathcal{R})$-basis for the  real $(n, m)$-module 
$\mathcal{R}^{(2)n, m}$.
For 
$$
x=\displaystyle\sum_{i=1}^n (\underbrace{x^{i1}+1^{\#}x^{i2}}_{x^i})e_i+\sum_{j=1}^m\,x^{n+j}e_{n+j}\in \mathcal{R}^{(2)n,m},
$$
$$
y=\displaystyle\sum_{i=1}^n (\underbrace{y^{i1}+1^{\#}y^{i2}}_{y^i})e_i+\sum_{j=1}^m\,y^{n+j}e_{n+j}\in \mathcal{R}^{(2)n,m},
$$
we define a real-valued function
$$<\, ,\,>: \mathcal{R}^{(2)n,m}\times \mathcal{R}^{(2)n,m}\to \mathcal{R}\quad\mbox{ and}\quad ||\,\,||: \mathcal{R}^{(2)n,m}\to \mathcal{R}$$ by
$$
<x, y>:=2\displaystyle\sum_{i=1}^n x^{i1}y^{i1}+\displaystyle\sum_{i=1}^n x^{i2}y^{i2}
+\displaystyle\sum_{j=1}^m x^{n+j}y^{n+j},
$$
where  $x^{i1}$, $x^{i2}$, $y^{i1}$, $y^{i2}$, $x^{n+j}$,  $y^{n+j}\in \mathcal{R}$ with
$1\le i\le n$ and $1\le j\le m$.  Then $<\, ,\,>$ is an inner product on $\mathcal{R}^{(2)n,m}$. Hence, 
$\mathcal{R}^{(2)n,m}$ is a metric space with the distance function $||x||:=\sqrt{<x, x>}$, where 
$x\in \mathcal{R}^{(2)n,m}$. 

\begin{definition}\label{def2}  Let $U$ be an open subset of $\mathcal{R}^{(2)n,m}$.  A function 
$f: U\to \mathcal{R}^{(2)s,t}$ is 
{\bf dual real differentiable} at $a\in U$ if there is a 
$\mathcal{R}^{(2)}$-module map $\lambda: \mathcal{R}^{(2)n,m}\to \mathcal{R}^{(2)s,t}$ such that
$\displaystyle\lim_{x\to a}\displaystyle\frac{||f(x)-f(a)-\lambda (x-a)||}{||x-a||}=0.$ 
\end{definition}

\medskip
The $\mathcal{R}^{(2)}$-module map $\lambda$ is denoted by $\mathcal{D}f(a)$ and 
called the {\bf $\mathcal{R}^{(2)}$-derivative} of $f$ at $a$. If $f$ is dual real differentiable at each point of the open subset $U$, then $f$ is said to be {\bf  dual real differentiable} on $U$. The familiar properties for the derivatives in Calculus hold for the  $\mathcal{R}^{(2)}$-derivative defined as a  $\mathcal{R}^{(2)}$-module map.
If $f: \mathcal{R}^{(2)n,m}\to \mathcal{R}^{(2)}=\mathcal{R}^{(2)1,0}$ is a function, then
$f(x)=(Re\,f)+1^{\#}(Ze\,f)$ for
$x=\displaystyle\sum_{i=1}^n (\underbrace{x^{i1}+1^{\#}x^{i2}}_{x^i})e_i+\sum_{j=1}^m\,x^{n+j}e_{n+j}$, where 
both $Re\,f$ and $Ze\,f$ are real-valued functions of $2n+m$ real variables:
$x^{11}, \dots, x^{n1}, x^{12}, \dots, x^{n2}, x^{n+1}, \dots, x^{n+m}$. We say that 
$f: \mathcal{R}^{(2)n,m}\to \mathcal{R}^{(2)}$ is {\bf  $\mathcal{R}^{(2)}$-smooth} if both $Re\,f$ and $Ze\,f$ are 
smooth. We say that a function $f: \mathcal{R}^{(2)n,m}\to  \mathcal{R}^{(2)s,t}$ is  {\bf  $\mathcal{R}^{(2)}$-smooth} if all of its component functions $f^1, \dots , f^s, f^{s+1}, \dots,  f^{s+t}$ are $\mathcal{R}^{(2)}$-smooth.

\medskip
We introduce the concept of zero divisor manifolds in the following

\begin{definition}\label{def2.2} A Hausdorff topological space $M$ is called a $(n, m)$-dimensional 
{\bf $\mathcal{R}^{(2)}$-smooth zero divisor manifold} if there exist open subsets $U_{\lambda}$ of $M$ with 
$\lambda\in \Lambda$, where  $\Lambda$ is a finite or countable set of indices, and  a map $\Phi_{\lambda}: U_{\lambda}\to \mathcal{R}^{(2)n,m}$  for every $\lambda\in \Lambda$  such that
\begin{description} 
\item[(i)] $M=\displaystyle \bigcup _{\lambda\in \Lambda} U_{\lambda}$;
\item[(ii)] The image $G_{\lambda}=\Phi_{\lambda}(U_{\lambda})$ is an open set in $ \mathcal{R}^{(2)n,m}$;
\item[(iii)] The map $\Phi_{\lambda}$ is one-to-one;
\item[(iv)]  For any two sets $ U_{\lambda}$, $ V_{\mu}$, $\lambda, \mu\in \Lambda$ the images 
$\Phi_{\lambda}(U_{\lambda}\cap V_{\mu}), \Phi_{\mu}(U_{\lambda}\cap V_{\mu})\subseteq  \mathcal{R}^{(2)n,m}$
are open and the transition map 
$\Phi_{\mu}\circ \Phi_{\lambda}^{-1}: \Phi_{\lambda}(U_{\lambda}\cap V_{\mu})\to \Phi_{\mu}(U_{\lambda}\cap V_{\mu})$ is $\mathcal{R}^{(2)}$-smooth.
\end{description}
The collection $\{( U_{\lambda}, \Phi_{\lambda})\, |\, \lambda\in \Lambda\}$ is called a {\bf $\mathcal{R}^{(2)}$-smooth atlas.}
\end{definition}

\medskip
As an example of $\mathcal{R}^{(2)}$-smooth zero divisor manifolds, we now construct the  real
projective  $\mathcal{R}^{(2)}$-module, which generalizes the ordinary real projective space.
Let  $\mathcal{R}^{(2)(n+1), m+1}$ be the real $(n+1, m+1)$-dimensional 
$\mathcal{R}^{(2)}$-module, and let $V$ and $W$ be defined by
$$V:=
\left\{(x^{02}1^{\#},  \dots ,  x^{n2}1^{\#}\,| \,x^{n+1}1^{\#},  \dots ,  x^{n+m+1}1^{\#})\left |
\begin{array}{c}\{x^{i2}\}_{i=0}^n\subseteq \mathcal{R}\\ \\\{x^{n+j+1}\}_{j=0}^m\subseteq \mathcal{R}\end{array}\right.\right\}$$
and
$$W:=
\Big\{(x^0,  \dots ,  x^n\,| \,\underbrace{0,  \dots ,  0}_{m+1})\, |\,
\{x^0,  \dots ,  x^n\}\subseteq \mathcal{R}^{(2)}\Big\}.$$
Let
$\overline{\mathcal{R}^{(2)(n+1), m+1}}:=\mathcal{R}^{(2)(n+1), m+1}\setminus (V\cup W)$ be a subset of $\mathcal{R}^{(2)(n+1), m+1}$. We define an equivalence relation $\sim$ on the subset
$\overline{\mathcal{R}^{(2)(n+1), m+1}}$ of $\mathcal{R}^{(2)(n+1), m+1}$ as follows:
\begin{eqnarray*}
&&\mbox{$x\sim y$ iff there exist $s\in \mathcal{R}^{(2)}$ with $Re\,s\ne 0$ and $t\in  \mathcal{R}$ with $t\ne 0$ such}\\
&&\mbox{that $y^i=sx^i$ for $0\le i\le n$ and $y^{n+1+j}=tx^{n+1+j}$ for $0\le j\le m$,}
\end{eqnarray*}
where 
\begin{eqnarray*}
&&x=(x^1,  \dots ,  x^n\,| \,x^{n+1}1^{\#},  \dots ,  x^{n+m}1^{\#})\in \overline{\mathcal{R}^{(2)(n+1), m+1}},\\
&&x=(y^1,  \dots ,  y^n\,| \,y^{n+1}1^{\#},  \dots ,  y^{n+m}1^{\#})\in \overline{\mathcal{R}^{(2)(n+1), m+1}}.
\end{eqnarray*}
The {\bf  real projective  $\mathcal{R}^{(2)}$-module} $\mathcal{R}^{(2)}\mathcal{P}^{n,m}$ is defined as the quotient space of $\overline{\mathcal{R}^{(2)(n+1), m+1}}$ by the equivalence relation above, i.e., 
$$
\mathcal{R}^{(2)}\mathcal{P}^{n,m}:=\overline{\mathcal{R}^{(2)(n+1), m+1}}/\sim
=\big\{[x]\,|\, x\in \mathcal{R}^{(2)}\mathcal{P}^{n,m}\big\},
$$
where $[x]$ denote the  equivalence class of $x\in \overline{\mathcal{R}^{(2)(n+1), m+1}}$.
For $0\le i\le n$ and $0\le j\le m$, we define
$U_{ij}:=\{[x]\,|\, \mbox{$Re\,x^i\ne 0$ and $x^{n+1+j}\ne 0$}\}$ and a map 
$\phi_{ij}: U_{ij}\to \mathcal{R}^{(2)n,m}$ by
\begin{eqnarray*}
&&\phi_{ij}\Big([x^0,  x^1,\dots ,  x^n\,| \,x^{n+1}1^{\#},  \dots ,  x^{n+m+1}1^{\#}]\Big)\\
&=&\Big(\frac{x^0}{x^i},  \dots ,  \widehat{\frac{x^i}{x^i}}, \dots, \frac{x^n}{x^i}\,\left | \, \frac{x^{n+1}}{x^{n+1+j}}1^{\#},  \dots ,   \widehat{\frac{x^{n+j+1}}{x^{n+1+j}}1^{\#}} ,\dots \frac{x^{n+m+1}}{x^{n+1+j}}1^{\#}\Big)\right..
\end{eqnarray*}
Then the collection $\{( U_{ij}, \phi_{ij})\, |\, \mbox{$0\le i\le n$ and $0\le j\le m$  }\}$ is  a 
 $\mathcal{R}^{(2)}$-smooth atlas, and $\mathcal{R}^{(2)}\mathcal{P}^{n,m}$ is a 
$(n,m)$-dimensional  $\mathcal{R}^{(2)}$-smooth zero divisor manifold.

\bigskip
\section{Symplectic $\mathcal{R}^{(2)}$-Modules}

The concept of symplectic $\mathcal{R}^{(2)}$-modules is based on the following generalization of linear symplectic form in symplectic geometry.

\medskip
\begin{definition}\label{def3.1} A {\bf  Hu-Liu symplectic form} on   a 
$\mathcal{R}^{(2)}$-module $V$ is a function $\omega: V\times V\to \mathcal{R}^{(2)}$
satisfying the following properties:
\begin{description}
\item[(i)] $\omega(x+y, z)=\omega(x, z)+\omega(y, z)$ for  $x, y, z\in V$;
\item[(ii)] $ \omega(ax, y)=a\omega(x, y)$ for $x, y\in V$ and $a\in \mathcal{R}^{(2)}$;
\item[(iii)] $\omega(x, y)=-\omega(y, x)$ for $x, y\in V$;
\item[(iv)] If $v\in V$ and $Re\,\omega(v, x)=0$ for all $x\in V$, then $v\in  Ker\,1^{\#}_V$;
\item[(v)] If $v\in  Ker\,1^{\#}_V$ and $\omega(v, x)=0$ for all $x\in V$, then 
$v\in  Im\,1^{\#}_V$.
\end{description}
\end{definition}

\medskip
\begin{definition}\label{def3.2} A pair $(V, \omega)$ is called a {\bf symplectic 
$\mathcal{R}^{(2)}$-module }  if $V$ is a $\mathcal{R}^{(2)}$-module and $\omega$ is a Hu-Liu symplectic form on   the $\mathcal{R}^{(2)}$-module $V$.
\end{definition}

\medskip
The basic property of  symplectic $\mathcal{R}^{(2)}$-modules is the following

\medskip
\begin{proposition}\label{pr3.1} If $(V, \omega)$ is a finite-dimensional  symplectic 
$\mathcal{R}^{(2)}$-module, then $\dim_{\mathcal{R}^{(2)}}(V)=(2n,\,2m)$ for some non-negative integers $n$ and $m$, and there exists a $(\mathcal{R}^{(2)}, \mathcal{R})$-basis
$(e_1, f_1, \dots, e_n, f_n \,\|\, e_{n+1}, f_{n+1}, \dots, e_{n+m}, f_{n+m})$ such that 
\begin{description}
\item[(i)] $\omega(e_i, f_i)=1 \quad\mbox{for $1\le i\le n$}$;
\item[(ii)]  $\omega(e_{n+j}, f_{n+j})=1^{\#} \quad\mbox{for $1\le j\le m$}$;
\item[(iii)] $\omega(e_i, e_j)=0 \quad\mbox{for $1\le i, j\le n+m$}$;
\item[(iv)] $\omega(f_i, f_j)=0 \quad\mbox{for $1\le i, j\le n+m$}$;
\item[(v)] $\omega(e_i, f_j)=0 \quad\mbox{for $1\le i, j\le n+m$ and $i\ne j$}$.
\end{description}
\end{proposition}

\end{document}